\documentclass[11pt,lettersize,reqno,fullpage]{article}
\usepackage{amsmath,amssymb,amsfonts,mathrsfs,verbatim,enumitem,pstricks,amsthm, graphicx}
\usepackage{authblk}
 \usepackage{amscd}
 \usepackage{epsf}
 \usepackage{epsfig}
 \usepackage{graphics}
 \usepackage{graphicx}
 \usepackage{geometry}
 \usepackage{tikz} 
\usepackage{authblk}
\usepackage[all]{xy}
 \usepackage{geometry}
 \usepackage{tikz} 
 \usepackage{amssymb}
 \usepackage{textcomp}
 \usepackage{gensymb}
 \usepackage{tkz-euclide}
 \usepackage{verbatim}
 \usepackage[english]{babel}
 \usepackage{tabu}
 \usepackage{amsmath}
 \usepackage{hyperref}
\usepackage{centernot, xr, accents}
\usepackage{hyperref}
\usepackage[bottom]{footmisc}
\usepackage{footnote}
\usepackage{amsmath,amssymb,amsfonts,mathrsfs,verbatim,enumitem,pstricks,amsthm}
\usepackage[all]{xy}
\usepackage{centernot, xr, accents}
\usepackage{hyperref}
\usepackage{amsmath,stackengine}
\usepackage{scalerel,stackengine}
\stackMath
\newcommand\reallywidehat[1]{%
\savestack{\tmpbox}{\stretchto{%
  \scaleto{%
    \scalerel*[\widthof{\ensuremath{#1}}]{\kern-.6pt\bigwedge\kern-.6pt}%
    {\rule[-\textheight/2]{1ex}{\textheight}}
  }{\textheight}%
}{0.5ex}}%
\stackon[1pt]{#1}{\tmpbox}%
}

\def \hf{\hspace*{0.5cm}}                      
\def\bge{\begin{equation}}                
\def\ede{\end{equation}}                
\def\bgd{\begin{displaymath}}         
\def\edd{\end{displaymath}}            
\def\bgee{\begin{equation*}}           
\def\edee{\end{equation*}}           

\def \ni{\noindent}

\def\lra{\longrightarrow}

\def\BA{\begin{eqnarray}}
\def\EA{\end{eqnarray}}
\def\BAA{\begin{eqnarray*}}
\def\EAA{\end{eqnarray*}}
\def\Bal{\begin{align*}}
\def\Eal{\end{align*}}

\def \C{\mathbb{C}}

\def \P{\mathbb{P}}
\def \A{\mathcal{A}}

\def \f{\frac}

\def \25node{A_5} 
\def \62node{A_6}

\def \D{\mathcal{D}}

\def \ov{\overline}

\def \q{\rho}

\def \U{\mathcal{U}}

\def \Z{\mathrm{A}} 
\def \Y{\mathrm{B}}

\def \sq{s}
\def \rq{r}

\def \hxt2{\hat{x}_{t_2}} 
\def \hyt2{\hat{y}_{t_2}}

\newcommand\groupequation[2][21pt]{%
  \setbox0=\hbox{$\displaystyle#2$}%
  \stackengine{0pt}{\copy0}{%
    \makebox[\linewidth]{\hfill$\left.\rule{0pt}{\ht0}\right\}$\kern#1}}
    {O}{c}{F}{T}{L}
}

\title{{\LARGE Counting singular curves with tangencies}}
\author{Anantadulal Paul}
\date{}
\usepackage{hyperref}
\hypersetup{
	colorlinks,
	citecolor=red,
	filecolor=black,
	linkcolor=blue,
	urlcolor=black
}
\theoremstyle{plain}
\newtheorem{thm}{Theorem}[section]

\newtheorem{lmm}[thm]{Lemma}

                                                                                                                                                                                                                                                                                                                                                                                                                                                                                                                                                                                                                                                                                                                                                                                                                                                                                                                                                                                                                                                                                                                                                                                                                                                                                                                                                                                                                                                                                                                                                                                                                                                                                                                                                                                                                                                                                                                                                                                                                                                                                                                                                                                                                                                                                                                                                                                                                                                                                                                                                                                                                                                                                                                                                                                                                                                                                                                                                                                                                                                                                                                                                                                                                                                                                                                                                                                                                                                                                                                                                                                                                                                                                                                                                                                                                                                                                                                                                                                                                                                                                                                                                                                                                                                  \theoremstyle{plain}
                                                                                                                                                                                                                                                                                                                                                                                                                                                                                                                                                                                                                                                                                                                                                                                                                                                                                                                                                                                                                                                                                                                                                                                                                                                                                                                                                                                                                                                                                                                                                                                                                                                                                                                                                                                                                                                                                                                                                                                                                                                                                                                                                                                                                                                                                                                                                                                                                                                                                                                                                                                                                                                                                                                                                                                                                                                                                                                                                                                                                                                                                                                                                                                                                                                                                                                                                                                                                                                                                                                                                                                                                                                                                                                                                                                                                                                                                                                                                                                                                                                                                                                                                                                                                                                  \usepackage{blindtext}
\newcommand{\address}

\newcommand{\email}

\newtheorem{obs}[thm]{Observation}

\newtheorem*{mthm}{MAIN THEOREM}

\theoremstyle{definition}
\newtheorem{rem}[thm]{Remark} 
\newtheorem{defn}[thm]{Definition}
 \numberwithin{equation}{section}

 \newcommand{\field}[1]{\ensuremath{\mathbb{#1}}}

 \newcommand{\QQ}{\field{Q}}

 \newcommand{\PP}{\field{P}}
 \newcommand{\DD}{\mathcal{D}}

 \def \A{\mathrm{A}}
 
 \def \D{\mathcal{D}}
 \def \lra{\longrightarrow} 
 \def \hf{\hspace*{0.5cm}}

 \def \U{U}
 \def \C{\mathbb{C}}
 \def \H{\mathcal{S}}
 \def \q{f}
 
 \def \hf{\hspace*{0.5cm}}                      
\def\bge{\begin{equation}}                
\def\ede{\end{equation}}                
\def\bgd{\begin{displaymath}}         
\def\edd{\end{displaymath}}            
\def\bgee{\begin{equation*}}           
\def\edee{\end{equation*}}   
\def \f{\frac}
\def \ni{\noindent}
\def \sq{s}
\def \rq{r} 
\def \ov{\overline}
\def \Z{\mathrm{A}} 
\def \Y{\mathrm{B}}
 \theoremstyle{definition}
\usepackage{hyperref} 
 
\begin{document}
\maketitle


\begin{abstract} 
We obtain  a recursive formula for the characteristic number of degree $d$ curves in $\mathbb{P}^2$ with prescribed singularities 
(of type $A_k$) 
that are tangent to a given line. The formula is in terms of the characteristic number of curves with exactly those singularities.  
Combined with the results of S.~Basu and R.~Mukherjee (\cite{R.M}, \cite{B.R} and \cite{BM8}), 
this gives us a complete formula for the characteristic number of curves with $\delta$-nodes 
and one singularity of type $A_k$, tangent to a given line, provided $\delta+k \leq 8$.
We use a topological method, namely the method of ``dynamic intersections'' (cf.~Chapter 11 in \cite{F}) 
to compute the degenerate contribution to the Euler class. Till codimension eight, we verify that our numbers are logically 
consistent with those computed earlier by Caporaso-Harris(\cite{CH}). 
We also make a non trivial low degree check to verify our formula for the number 
of cuspidal cubics tangent to a given line, using a result of Kazaryan (\cite{Kaz}).
\end{abstract}
\tableofcontents

 	\section{Introduction}
 	A fundamental problem in enumerative geometry is to count curves with prescribed singularities. 
This question has been studied for a very long time starting with Zeuthen (\cite{Zu}) more than a hundred years ago. 
It has been studied extensively in the last thirty years from various perspectives by numerous mathematicians including amongst others,  
Z.~Ran (\cite{Ran1}, \cite{ZRan}), I.~Vainsencher (\cite{Van}), L.~Caporaso and J.~Harris (\cite{CH}), M.~Kazaryan (\cite{Kaz}), S.~Kleiman and R.~Piene (\cite{KP1}), 
D.~Kerner (\cite{Ker1} and \cite{Ker2}), 
F.~Block (\cite{FB}),
Y.~J. Tzeng and J.~Li (\cite{Tz}, \cite{Tzeng_Li}), M.~Kool, V.~Shende and R.P.~Thomas (\cite{KST}), S.~Fomin and G.~Mikhalkin (\cite{FoMi}), 
G.~Berczi (\cite{Berczi}) and S.~Basu and R.~Mukherjee (\cite{R.M}, \cite{B.R} and \cite{BM8}).\\
\hf \hf A closely related question is to enumerate curves with prescribed singularities that are tangent to a given line.  
This question also has a long history that can be traced back to Zeuthen. As early as   
$1848$, Zeuthen computed the characteristic 
number of rational quartics in $\mathbb{P}^2$ 
tangent to a given line. \\ 
\hf \hf In the last thirty years, an extensive amount of work has been done in enumerating curves that are tangent to a given 
line when the \textit{prescribed singularities are nodes}. These include among others the results of Z.~Ran (\cite{Ran}), 
I.~Vainsencher(\cite{Van}), Caporaso-Harris (\cite{CH}), R.~Vakil(\cite{rv1}, \cite{rv2}), A.~Gathmann (\cite{AnGa1} and  \cite{AnGa2}) and C.~Cadman and L.~Chen (\cite{CadLi}).\\  
\hf \hf Very recently, using methods of algebraic cobordism, Y.~J. Tzeng has shown (\cite{Tzeng_Tang}) that a universal formula exists for the characteristic 
number of curves in a linear system, that are tangent to a given line and that have prescribed singularities (more degenerate than nodes). \\
\hf  \hf   We now mention a result that we are aware of concerned with the tangency question in other spaces. In \cite{Y.Cooper}, Y. Cooper and R. Pandharipande study the Severi problem involving single tangency condition via the matrix elements in Fock space.\\
\hf  \hf In this paper, we obtain a recursive formula  for the characteristic number of curves that are tangent to a given line and that 
have prescribed singularities (of type $A_k$). Furthermore, till codimension eight we can obtain explicit formulas. The method we use is the method of 
dynamic intersection theory, similar to what is used in \cite{Zin}, \cite{R.M}, \cite{B.R} and \cite{BM8}.
Before stating the main result of the paper, let us make a couple of definitions: 
\begin{defn}
\label{singularity_defn}
Let $f:\mathbb{P}^2  \longrightarrow \mathcal{O}(d)$ be a holomorphic section. 
A point $q \in f^{-1}(0)$ \textsf{is of singularity type} $A_k$
if there exists a coordinate system
$(x,y) :(U,q) \lra (\C^2,0)$ such that $f^{-1}(0) \cap \U$ is given by 
\begin{align*}
y^2 + x^{k+1}   &=0.   
\end{align*}
\end{defn}

In more common terminology, $q$ is a {\it smooth} point of $f^{-1}(0)$ if 
it is a singularity of type $\A_0$, a {\it simple node} (or just node) if its singularity type is $\A_1$, 
a {\it cusp} if its type is $\A_2$ and a {\it tacnode} if its type is $\A_3$. 
We will frequently use the phrase ``a singularity of codimension $k$''.  
This refers to the number of conditions having that singularity imposes on the space of curves. 
More precisely, it is the expected codimension of the equisingular strata.  
Hence, an $A_k$ singularity is a singularity of codimension $k$.\\  
\hf \hf Next, given a non negative integer $k$ and 
positive integers $\delta_1, \ldots, \delta_k$, let us define  
$N_{d}(A_1^{\delta_1} \ldots A_{k}^{\delta_k})$ to be the number of degree $d$-curves in $\mathbb{P}^2$, passing through 
$\frac{d(d+3)}{2}-  (\delta_1 + 2 \delta_2 + \ldots  + k \delta_{k})$ generic points having $\delta_i$ ordered singularities 
of type $A_i$. \\ 
\hf \hf Similarly, we define 
$N_{d}(A_1^{\delta_1} \ldots A_{k}^{\delta_k}; L_{A_i})$ to be the number of 
degree $d$-curves in $\mathbb{P}^2$, passing through 
$\frac{d(d+3)}{2}-  (1+\delta_1 + 2 \delta_2 + \ldots  + k \delta_{k})$ generic points, having $\delta_j$ ordered singularities 
of type $A_j$ (when $i \neq j$), $\delta_i-1$ singularities of type $A_i$ and another singularity of type $A_i$ lying on a given line. \\
\hf \hf The main result of this paper is as follows: 
\begin{mthm}
\label{main_thm}
Let $k$ be a non negative integer and 
$\delta_1, \delta_2, \ldots, \delta_k$ a collection of positive integers. 
Define 
\[ \delta_d:= \frac{d(d+3)}{2} \qquad \textnormal{and} \qquad w_{d}:= \delta_d-  (1+\delta_1 + 2 \delta_2 + \ldots  + k \delta_{k}). \]
Let $N^{\mathrm{T}}_d(A_1^{\delta_1} \ldots A_k^{\delta_{k}})$ 
denote the number 
of degree $d$-curves in $\mathbb{P}^2$, passing through 
$w_d$ generic points, 
having $\delta_i$ (ordered) singularities of type $A_i$ 
(for all $i$ from $1$ to $k$) that is tangent to a given line. 
Then, 
\begin{align}
N^{\mathrm{T}}_d(A_1^{\delta_1} \ldots A_k^{\delta_{k}})& =
 2(d-1){N}_d(A_1^{\delta_1} \ldots A_k^{\delta_{k}}) - {\sum_{i=1}^{i=k}{{\delta_i }
 (i+1)N_d(A_1^{\delta_1} \ldots A_k^{\delta_{k}}; L_{A_i}}}), \label{main_formula}
\end{align} 
for all $d \geq d_{\textnormal{min}}$,  where 
\[ d_{\textnormal{min}}:= k+(2\delta_1+ \delta_2+ \ldots + \delta_k).\]
\end{mthm}

\begin{rem}
We note that the numbers $N_d(A_1^{\delta} A_k)$ are directly given in the papers of 
S.~Basu and R.~Mukherjee (\cite{R.M}, \cite{B.R} and \cite{BM8}) when $\delta+k \leq 8$. 
The results of those papers can be used to compute 
$N_d(A_1^{\delta} A_k; L_{A_i})$ when $\delta+k \leq 8$ with no further effort (since they obtain
 equality on the level of cycles). Hence,  
using these numbers and using equation (\ref{main_formula})
we can obtain a complete formula for $N_d^{\mathrm{T}}(A_1^{\delta} A_k)$ when $\delta+k \leq 8$. 
The formulas for 
$N_d^{\mathrm{T}}(A_k)_{1 \leq k \leq 8}$, $N_d^{\mathrm{T}}(A_1 A_k)_{1 \leq k \leq 7}$ and $N_d^{\mathrm{T}}(A_1^{\delta})_{1 \leq \delta \leq 8}$
are listed explicitly in section \ref{ef}.
\end{rem}

\begin{rem}
Next, we note that in \cite{Kaz}, M.~Kazaryan computes all the characteristic number of curves with upto 
seven singularities. We believe he obtains equality on the level of cycles; hence we believe in principle his method can be used 
to compute the characteristic number of curves with singularities, where one of the singularity is required to 
lie on a line (till codimension seven).\footnote{However, we are not completely certain about this point.} 
Hence, using equation (\ref{main_formula}), we can in principle obtain 
a formula for $N_d^{\mathrm{T}}(A_1^{\delta_1} A_2^{\delta_2} \ldots A_k^{\delta_k})$ 
provided the total codimension is seven. 
\end{rem}

\begin{rem}
\label{delta_0}
When $k=0$, we will abbreviate $N_d(A_1^{\delta_1} \ldots A_k^{\delta_{k}})$ as $N_d$ 
and we will abbreviate 
$N^{\mathrm{T}}_d(A_1^{\delta_1} \ldots A_k^{\delta_{k}})$ as $N^{\mathrm{T}}_d$. 
We note that $N_d$ is the 
number 
of degree $d$ curves in $\mathbb{P}^2$ passing through $\delta_d$ generic points; hence 
$N_d =1$. Similarly, $N^{\mathrm{T}}_d$
is the number 
of degree $d$ curves in $\mathbb{P}^2$ passing through $\delta_d-1$ generic points that is tangent to a given line.  
Hence, the $k=0$ case of equation \eqref{main_formula} implies 
\[ N_d^{\mathrm{T}} = 2(d-1).\]
\end{rem}

\begin{rem}
The bound $d \geq d_{\textnormal{min}}$ is imposed to ensure we get transversality of certain sections. However, this bound is not necessarily sharp; 
the bound is sufficient to get transversality, but it is not always necessary. 
\end{rem}

\section{Overview of the method}

We now give an overview of the method we use. Our starting point will be the following classical fact from Differential Topology: 
\begin{thm}
\label{Fact}
Let $V \longrightarrow M$ be an oriented vector bundle over a compact, oriented manifold M 
and $s:M \longrightarrow V$ a section that is transverse to zero. 
 
If the rank of $V$ is equal to the dimension of $M$,
then the signed cardinality of $s^{-1}(0)$ is the Euler class of $V$, evaluated on the fundamental class of M, i.e.,\\
\hspace*{2cm} \[\mid \pm s^{-1}(0) \mid = \left\langle e(V), [M]\right\rangle\]
\end{thm}

\begin{rem}
We will express the tangency condition as the vanishing of a section of an appropriate vector bundle. However, the corresponding 
Euler class involves a degenerate contribution. The central aspect of this paper is how we compute the degenerate contribution 
to the Euler class. 
We use the method of ``dynamic intersections'' (cf.~Chapter 11 in \cite{F}) 
to compute the degenerate contribution to the Euler class.
 
\end{rem}

\section{Proof of Main Theorem}

Let us denote $\D$ to be the space of non-zero homogeneous degree $d$-polynomials in three variables upto scaling, i.e., 
\[ \D := \mathbb{P}(H^0(\mathbb{P}^2, \mathcal{O}(d))) \approx \mathbb{P}^{\delta_d}.\]

Hence, $\D$ can be identified with the space of degree $d$ curves in $\mathbb{P}^2$ (not necessarily irreducible).

Let 
\begin{align*}
\gamma_{\D} &\longrightarrow \D \qquad \textnormal{and} \qquad \gamma_{\mathbb{P}^2} \longrightarrow \mathbb{P}^2 
\end{align*}
be the tautological line bundles over $\D$ and $\mathbb{P}^2$ respectively.  \\

\hf \hf We will now prove our main theorem (i.e., we will prove \eqref{main_formula}). 

Given non negative integer $k$ and positive integers 
$\delta_1, \ldots \delta_{k}$, let us define 

\begin{align}
M &:= \D \times (\mathbb{P}^{2})^{\delta_1} \times \ldots \times (\mathbb{P}^{2})^{\delta_k} \qquad \textnormal{and} \nonumber \\
\mathcal{S} &:= \{ ([f], q_1^1,\ldots q^1_{\delta_1};  \ldots; q^k_{1}, \ldots q^k_{\delta_k}) \in M: ~~
\textnormal{$f$ has an $A_i$ singularity at $q^i_{\eta}$}, ~~q^{i}_{\eta} ~~\textnormal{are all distinct}\}. \label{s_manifold}
\end{align}
We will show shortly that if $d \geq d_{\textnormal{min}}$, then $\H$ is a complex sub manifold of $M$, of dimension 
$w_d+1$. Let us now  make the following abbreviation: 
\begin{align*}
\overline{q}&:=  (q_1^1,\ldots q^1_{\delta_1};  \ldots; q^k_{1}, \ldots q^k_{\delta_k}) \in (\mathbb{P}^{2})^{\delta_1} \times \ldots \times (\mathbb{P}^{2})^{\delta_k}.
\end{align*}

We now define the following two sections of line bundles over $\overline{\H}\times L$:
 
\begin{align*}
\psi_{\textnormal{ev}}: \overline{\H}\times L & \longrightarrow \mathbb{L}_{\textnormal{ev}}:= \gamma_{\D}^*\otimes \gamma_{L}^{*d}, \qquad \textnormal{given by} \qquad 
\{\psi_{\textnormal{ev}} ([f], \overline{q}, p)\}(f) := f(p) \qquad \textnormal{and} \\ 
\psi_{\textnormal{T}}: \psi_{\textnormal{ev}}^{-1}(0) & \longrightarrow \mathbb{L}_{\textnormal{T}}:= \gamma_{\D}^*\otimes T^*L \otimes \gamma_{L}^{*d}, 
\qquad \textnormal{given by} \qquad 
\{\psi_{\textnormal{T}} ([f], \overline{q}, p)\}(f \otimes v) := \nabla f|_p (v). 
\end{align*}
Here $\gamma_{L}$ denotes the tautological line bundle over $L$ 
(which is the same as the restriction of the tautological line bundle $\gamma_{\mathbb{P}^2}$ to $L$). \\
\hf \hf Let us now define 
\begin{align*}
\mathcal{B}^i_{\eta}&:= \{([f], \overline{q}, p) \in \overline{\H}\times L: q^i_{\eta} = p\} \qquad \textnormal{and} \qquad \mathcal{B}:= \bigcup \mathcal{B}^i_{\eta}. 
\end{align*}
We claim that restricted to $\H\times L-\mathcal{B}$, the sections  
 
$\psi_{\textnormal{ev}}$ and $\psi_{\textnormal{T}}$ 
are transverse to zero. We will prove that claim shortly. \\ 
 
\hf \hf Next, let $\mu$ be the subspace of curves in $\D$ that pass through  $w_d$
generic points and let 
\[\pi_{\mathcal{D}}:M \times L \longrightarrow \D\]
be the projection map. 

Since the points are generic, the sub manifold $\pi_{\D}^{-1}(\mu)$ will intersect $\H \times L$ transversally (inside 
$M \times L$).\\ 

\hf \hf Next, we note that if $f$ is tangent to $L$ at $p$, then  
\begin{align}
f(p) &=0 \qquad \textnormal{and} \qquad \nabla f|_{p} (v) =0 \qquad \forall ~v \in T_pL. \label{psi_T_1}
\end{align}
In other words, 
\begin{align}
\psi_{\textnormal{ev}} ([f], \overline{q}, p) &=0 \qquad \textnormal{and} \qquad \psi_{\textnormal{T}} ([f], \overline{q}, p) =0.  \label{psi_T_2}
\end{align}

However, equation \eqref{psi_T_2} is also satisfied on $\mathcal{B}$ 
(i.e., when one of the singular points $q^{i}_{\eta}$, happens to lie on the line $L$, i.e., one of the points $q^{i}_{\eta}$ becomes 
equal to the tangency point $p$).  Hence, our desired number 
$N^{\mathrm{T}}_d(A_1^{\delta_1} \ldots A_k^{\delta_{k}})$ is the number of solutions to  
\begin{align}
\psi_{\textnormal{ev}} ([f], \overline{q}, p) &=0, \qquad \psi_{\textnormal{T}} ([f], \overline{q}, p) =0, \qquad 
([f], \overline{q}, p) \in \Big(\H \times L -\mathcal{B}\Big)\cap \Big(\pi_{\D}^{-1}\mu\Big). \label{num_of_soln_open}
\end{align}
We note that since the points are generic, 
\begin{align*}
\Big(\H \times L -\mathcal{B}\Big)\cap \Big(\pi_{\D}^{-1}\mu\Big) & =  \Big(\overline{\H} \times L -\mathcal{B}\Big)\cap \Big(\pi_{\D}^{-1}\mu\Big)
\end{align*}

Hence, we conclude that 
\begin{align}
\langle e(\mathbb{L}_{\textnormal{ev}})e(\mathbb{L}_{\textnormal{T}}), 
~[\overline{\H}\times L] \cap [\pi_{\D}^{-1}\mu] \rangle &=   
N^{\mathrm{T}}_d(A_1^{\delta_1} \ldots A_k^{\delta_{k}}) + \mathcal{C}_{\mathcal{B}_{\mu}}, \label{Euler_eqn}
\end{align}
where $\mathcal{C}_{\mathcal{B}_{\mu}}$ is the contribution of the section from the boundary $\mathcal{B}\cap \big(\pi_{\D}^{-1}\mu\big)$.
We note that $\cap$ denotes intersection inside the space $M\times L$. \\
\hf \hf Next, we note that the left hand side of equation \eqref{Euler_eqn} is given by 
\begin{align}
\langle e(\mathbb{L}_{\textnormal{ev}})e(\mathbb{L}_{\textnormal{T}}), 
~[\overline{\H}\times L] \cap [\pi_{\D}^{-1}\mu] \rangle & = 2(d-1){N}_{d}(A_1^{\delta_1} \ldots A_k^{\delta_{k}}). \label{euler}
\end{align}
We will now compute the quantity  
$\mathcal{C}_{\mathcal{B}_{\mu}}$. Let us first analyze the set $\mathcal{B}\cap \big(\pi_{\D}^{-1}\mu\big)$. 
This is the union 
of the sets 
\begin{align*}
\mathcal{B}^i_{\eta}\cap \big(\pi_{\D}^{-1}\mu\big). 
\end{align*}
We now note that  
$\mathcal{B}^i_{\eta}\cap \big(\pi_{\D}^{-1}\mu\big)$ 
is the set of all degree $d$ curves passing through the  
$w_d$ generic points, having $\delta_j$ (ordered) singularities of type $A_j$ (for all $j$ from $1$ to $k$) 
and where the $(q^{i}_{\eta})^{\textnormal{th}}$ singular point lies on a line. 
Note that the $(q^{i}_{\eta})^{\textnormal{th}}$ singular point corresponds to a singularity of type $A_i$. 

\begin{align*}
|\mathcal{B}^i_{\eta}\cap \big(\pi_{\D}^{-1}\mu\big)| &= N_d(A_1^{\delta_1} \ldots A_k^{\delta_{k}}; L_{A_i}).
\end{align*}
We claim that each point of $\mathcal{B}^i_{\eta}\cap \big(\pi_{\D}^{-1}\mu\big)$ vanishes with a multiplicity of $(i+1)$. 
Hence, the total contribution from the set $\mathcal{B}\cap \big(\pi_{\D}^{-1}\mu\big)$ to the Euler class is given by 
\begin{align}
\mathcal{C}_{\mathcal{B}_{\mu}} & = {\sum_{i=1}^{i=k}{{\delta_i }
 (i+1)N_d(A_1^{\delta_1} \ldots A_k^{\delta_{k}}; L_{A_i}}}). \label{bdry}
\end{align}
Equations \eqref{Euler_eqn}, \eqref{euler} and \eqref{bdry} give us equation \eqref{main_formula}. \\ 
\hf \hf We will now prove the claims that we have made regarding transversality and multiplicity. 

\section{Transversality and Multiplicity}

We will start by recalling a few facts about $A_k$ singularities that are proved in \cite{R.M}, section 3. 

Let $\mathcal{U}$ be a 
neighbourhood of the origin in $\C^2$ and 
$f:\mathcal{U} \longrightarrow \mathbb{C}$ 
a holomorphic function. 
Let $i,j$ be non-negative integers. We define
\bgd
\q_{ij}:=\frac{\partial^{i+j} \q}{\partial^i x\partial^j y}\bigg|_{(x,y)=(0,0)}\,.
\edd

The procedure to obtain $A^f_k$ is given in the proof of the following Proposition. 
We will now state a necessary and sufficient criteria for a curve to 
have a specific singularity of type $A_{k\geq 1}$.
\begin{lmm}
\label{ift}
Let $\q =\q(x,y)$ be a holomorphic function defined on a neighbourhood
of the origin in $\C^{2}$ such that $\q_{00} =0$ and $\nabla \q |_{(0,0)} \neq 0$. 
Then the origin is a smooth point of the curve. 

\end{lmm}

\begin{lmm}
\label{ml}
Let $\q = \q(x,y)$ be a holomorphic function defined on a neighbourhood
of the origin in $\C^2$ such that $\q_{00} =0$, $\nabla \q |_{(0,0)} =0$ and
$\nabla^2 \q|_{(0,0)}$ is non-degenerate.  
Then the curve has a singularity of type $A_1$ at the origin. 
\end{lmm}

\begin{rem}
Lemma \ref{ift} is also known as the \textit{Implicit Function Theorem} and Lemma \ref{ml} 
is also known as the \textit{Morse Lemma}. 
\end{rem}

\ni We now state a necessary and sufficient condition for a curve to have an $A_{k \geq 2}$ singularity. 
This can be thought of as a continuation of Lemma \ref{ml}. 

\begin{lmm}\label{fstr_prp}
Let $\q =\q(\rq, \sq)$ be a holomorphic function defined on a neighbourhood of the origin in $\C^{2}$ such that 
$\q_{00} =0$, $\nabla \q|_{(0,0)}=0$ and there exists a non-zero vector $\eta=(v_1,v_2)$ such that at 
the origin $ \nabla^2 f (\eta, \cdot)=0$, i.e., the Hessian is degenerate.
Let $x := v_1 \rq + v_2 \sq, y := -\ov{v}_2 \rq + \ov{v}_1 \sq $ and $\q_{ij}$ be the partial derivatives with 
respect to the new variables $x$ and $y$. 
Then, the curve $\q^{-1}(0)$ has a singularity of type $A_k$ at the origin  
if $\q_{02} \neq 0$ and 
the directional derivatives 
$A^{\q}_i$ defined in \eqref{Ak_sections} are zero for all $i \leq k$ and $A^{\q}_{k+1} \neq 0$. 
\end{lmm}
  
{\bf Proof}: The result follows from the following observation.
\begin{obs}
Let $\q =\q(\rq, \sq)$ be a holomorphic function defined on a neighbourhood of the origin in $\C^{2}$ such that 
$\q(0,0), ~\nabla \q|_{(0,0)}=0$ and there exists a non-zero vector $\eta = (v_1, v_2)$ such that at the 
origin $ \nabla^2 f (v, \cdot)=0$, i.e., the Hessian is degenerate.
Let $x := v_1 \rq + v_2 \sq, y := -\ov{v}_2 \rq + \ov{v}_1 \sq $ and $\q_{ij}$ be the partial derivatives 
with respect to the new variables $x$ and $y$. 
If $\q_{02} \neq 0$, there exists a coordinate chart $(u, v)$ centered around the 
origin in $\C^2$ such that 
\bge
\label{ak1}
\q = \left\{\begin{array}{rl}
v^2, & \textup{or}\\
v^2 + u^{k+1}, & \textup{for some $k\geq 2$.}
\end{array}\right.
\ede
\end{obs}
In terms of the new coordinates we have $\q_{00}= \q_{10}= \q_{01}= \q_{20}= \q_{11} =0$ and $\q_{02} \neq 0.$ 
Here $\partial_x + 0 \partial_y = (1,0)$ is the distinguished direction along which the Hessian is degenerate. \\

{\bf Proof of observation}: Let the Taylor expansion of $\q$ in the new coordinates be given by 
\bgd
\q(x,y) = \Z_0(x) + \Z_1(x)y + \Z_2(x) y^2 + \ldots.
\edd
By our assumption on $\q$, $\Z_2(0) \neq 0.$ We claim that there exists a holomorphic function $\Y(x)$ such that 
after we make a change of coordinates $y = y_1 + \Y(x)$, the function $\q$ is given by 
\bgd
\q = \hat{\Z}_0(x) + \hat{\Z}_2(x) y_1^2 + \hat{\Z}_3(x) y_1^3 + \ldots 
\edd
for some $\hat \Z_k(x)$ (i.e., $\hat{\Z}_1(x) \equiv 0$). To see this, 
we note that this is possible if $\Y(x)$ satisfies the identity
\begin{align}
\Z_1(x) + 2 \Z_2(x) \Y + 3 \Z_3(x) \Y^2 + \ldots \equiv 0.  \label{psconvgg}
\end{align}
Since ~$\Z_2(0) \neq 0$, $\Y(x)$ exists by the Implicit Function Theorem
\footnote{Moreover it is unique if we require $\Y(0) =0$.}.
Therefore, we can compute $\Y(x)$ as a power series using \eqref{psconvgg} and then
compute $\hat \Z_{0}(x)$. Hence, 
\begin{align}
\label{Ak_node_conditionn}
\q &= v^2 + \frac{A^{\q}_3}{3!}x^3 + 
\frac{A^{\q}_4}{4!} x^4 + \ldots,  
~~\textnormal{where} ~~~ v = \sqrt{(\hat{\Z}_2 + \hat{\Z}_3 y_1 + \ldots)} y_1, 
\end{align}
satisfies \eqref{ak1}. \qed \\
\hf\hf Following the above procedure we find $A_i^f$ for any $i$. For example,  
\begin{align}
\label{Ak_sections}
A^{\q}_3&= \q_{30}\,,\qquad
A^{\q}_4 = \q_{40}-\frac{3 \q_{21}^2}{\q_{02}}\,, \qquad
A^{\q}_5= \q_{50} -\frac{10 \q_{21} \q_{31}}{\q_{02}} + 
\frac{15 \q_{12} \q_{21}^2}{\q_{02}^2}, \ldots 
%
\end{align} 
and so on. We are now ready to prove the claim that the space of curves with prescribed singularities is a smooth manifold  
of the expected dimension, provided $d$ is sufficiently large. 

\begin{lmm}
Let $M$ and $\H$ be as defined in equation \eqref{s_manifold}. 
If $d \geq  d_{\textnormal{min}}$, then $\H$ is a complex sub manifold of $M$, of dimension 
$w_d+1$.
\end{lmm}
\noindent \textbf{Proof:} We will prove this statement by considering an affine chart. Hence, let us consider the 
vector space $\mathcal{F}_d \approx \mathbb{C}^{\frac{d(d+3)}{2}+1}$ of polynomials in two variables of degree at most $d$. 
Let us denote $p_i:= (x_i, y_i) \in \mathbb{C}^2$
and define 
\begin{align*}
\H_{\textnormal{affine}} &:= \{(f, p_1, p_2, \ldots, p_{\delta}) \in \mathcal{F}_d \times \Big((\mathbb{C}^2)^{\delta}-\Delta\Big): 
\textnormal{$f$ has an $A_{k_i}$-singularity at $p_i$, ~$p_i$ all distinct}\}. 
\end{align*}
Here $\Delta$ denotes the fat diagonal of $(\mathbb{C}^2)^{\delta}$ (i.e., if any two points are equal,  they belong to the fat diagonal). 
We will show that $\H_{\textnormal{affine}}$ is a smooth complex sub manifold of 
$\mathcal{F}_d \times \Big((\mathbb{C}^2)^{\delta}-\Delta\Big)$ of codimension $c_d$. 
In order to do that, we will describe $\H_{\textnormal{affine}}$ locally as the zero set of certain holomorphic functions.\\ 
\hf \hf Let us suppose that 
\[ (f, \overline{p}):= (f, p_1, p_2, \ldots, p_{\delta}) \in \mathcal{S}_{\textnormal{affine}}\]  
Suppose $f$ has an $A_{k}$ singularity at $p_1:= (x_1, y_1)$, then we can use Lemma $\ref{fstr_prp}$ to see that there exist a
sufficiently small open sets $U_{p_{1}} \subset \mathcal{F}_d $ around $f$ and $V_{p_{1}} \subset \mathbb{C}^2$ around $p_1$
 such that on $U_{p_{1}} \times V_{p_{1}}$(possibly after making a linear change of coordinates) $f_{yy}(x_1, y_1)$, the 
 second partial derivative of $f$ with respect to $y$, evaluated at $(x_1, y_1)$ is non zero. 
Let us now define 
\[ \hat{x}:= x - \frac{f_{xy}(x_1,y_1)}{f_{yy}(x_1,y_1)} y \qquad \textnormal{and} \qquad \hat{y}:= y. \] 
We note that $\hat{x}$ is well defined, since $f_{yy}(x_1, y_1) \neq 0$. 
We will now define $A^{f(x_1,y_1)}_k$ to be the expressions obtained in \eqref{Ak_sections}, where we replace $f_{ij}$ with the 
$(i,j)^{\textnormal{th}}$ partial derivative of $f$ with respect to $\hat{x}$ and $\hat{y}$, evaluated at $(x_1,y_1)$. 
As an example, 
\begin{align*}
A^{f(x_1,y_1)}_3 & = \Big( \partial_{x} - \frac{f_{xy}(x_1,y_1)}{f_{yy}(x_1,y_1)} \partial_{y}\Big)^3 f \\
           & = \Big( f_{xxx} - 3 \frac{f_{xy}}{f_{yy}} f_{xyy} + 3 \Big(\frac{f_{xy}}{f_{yy}}\Big)^2 f_{xyy} + 
            \Big(\frac{f_{xy}}{f_{xy}}\Big)^3 f_{yyy}\Big)\Big|_{(x_1,y_1)}.
\end{align*}
Since $f$ has an $A_{k_i}$ singularity at $p_i:= (x_i, y_i)$, all are distinct points so we can assume (possibly after a
 linear change of co-ordinates) that  $f_{yy}(x_i, y_i)$, the second partial derivative of $f$ with respect to $y$,
   evaluated at $(x_i, y_i)$ is non zero. Then repeated use of Lemma \ref{fstr_prp} will give us sufficiently small open 
   neighbourhoods $U:= \cap_{i} U_{p_{i}} \subset \mathcal{F}_d$ and $V:= \prod_{i} V_{p_{i}} \subset \Big((\mathbb{C}^2)^{\delta}-\Delta\Big)$
   so that we can define $A_{k_i}^{f(x_i, y_i)}$ to be the expressions as obtained in \eqref{fstr_prp} for each $i$.\\
Next, let $\mathcal{U}:= U \times V$ be a sufficiently small open neighbourhood of $(f, \overline{p})$ in $\mathcal{F}_d \times \Big((\mathbb{C}^2)^{\delta}-\Delta\Big)$. 
Let us define the 
function $\Phi: \mathcal{U} \longrightarrow \mathbb{C}^{c_d}$, given by 
\begin{align*}
\Phi(f, \overline{p})& := \Big(f(x_1, y_1), f_{x}(x_1, y_1), f_{y}(x_1, y_1), A^{f(x_1, y_1)}_{2}, \ldots, A^{f(x_1, y_1)}_{k_1};  \\ 
                     &         \qquad f(x_2, y_2), f_{x}(x_2, y_2), f_{y}(x_2, y_2), A^{f(x_2, y_2)}_{2}, \ldots, A^{f(x_2, y_2)}_{k_2};  \ldots, ; \\ 
                     &         \qquad f(x_{\delta}, y_{\delta}), f_{x}(x_{\delta}, y_{\delta}), f_{y}(x_{\delta}, y_{\delta}), 
                               A^{f(x_{\delta}, y_{\delta})}_{2}, \ldots, A^{f(x_{\delta}, y_{\delta})}_{k_{\delta}} \Big). 
\end{align*}
We claim that $\mathbf{0}$ is a regular value of $\Phi$. If we can show that, then our claim is proved. \\ 
\hf \hf To prove the claim, we will construct curves. 
Since the points $p_i$ are all distinct, we will show that for different possibilities of points we can produce 

 curves $\eta_i \in \mathcal{F}_d$ be such that 
\begin{align*}
\eta_i(x_j, y_j) &= \delta_{i,j}. 
\end{align*}
\begin{rem}
There are plenty of ways one can construct such curves $\eta_i $. In practice 
it is enough to construct curves $\eta_i$ such that $\eta_i(x_j, y_j) \neq 0$.
\end{rem}
We can easily construct such an $\eta_i$ by taking product of all the $(x-x_j)$, except $(x-x_i)$ combined with $(y-y_j)$, i.e., for 
$n$ distinct points
\begin{align*}
\eta_i & := 
                (z_{r_{1}}-z_1)\cdots (z_{r_{i-1}}-z_{i-1}) \reallywidehat{(z_{r_{i}}- z_i)} (z_{r_{i+1}}-z_{i+1}) \cdots (z_{r_{n}}-z_n)           
\end{align*}

 where \[(z_{r_{s}}, z_{s})=
                \begin{cases}
                (x, x_s) \hspace*{.5cm}  if \hspace*{.5cm} x_i \neq x_s \\
                (y, y_s) \hspace*{.5cm}  if \hspace*{.5cm} y_i \neq y_s
               \end{cases}\]
Let us consider the point $p_i:= (x_i, y_i)$. The curve $f$ has an $A_{k_i}$ singularity at $p_i$. 
As an example, if $f$ has at least $A_1$ singularity at $p_i$ then there are sufficiently small neighbourhoods 
around each $p_i$ where $f(p_i), f_x(p_i), f_y(p_i)$ vanishes. So in this
situation, if we simply construct curves as follows:
 \begin{align*}
\gamma^{i}_{00}(t) &:= f+ t \eta^{2}_i, ~~\gamma^{i}_{10}(t):=f + t (x-x_i) \eta^{2}_i, ~~\gamma^{i}_{01}(t):=f + t (y-y_i) \eta^{2}_i.
\end{align*}
So the above construction enables us
\begin{align*}
\begin{split}
\{d\Phi|_{(f, \overline{p})}\} (\gamma^{i^{\prime}}_{{0}, {0}}(0)) = (0, \cdots, \underbrace{(*,0,0)}_\text{i th position}, \cdots, 0)\\
\{d\Phi|_{(f, \overline{p})}\} (\gamma^{i^{\prime}}_{{1}, {0}}(0)) = (0, \cdots, \underbrace{(0,*,0)}_\text{i th position}, \cdots, 0)\\
\{d\Phi|_{(f, \overline{p})}\} (\gamma^{i^{\prime}}_{{0}, {1}}(0)) = (0, \cdots, \underbrace{(0,0,*)}_\text{i th position}, \cdots, 0)\\
\end{split}
\end{align*}
then the above computation implies that $\mathbf{0}$ is a regular value as claimed.\\
\hf Next, note that if $f$ has a singularity at least as degenerate as cusp at some point assuming that there is already $A_1$ singularity present at that point, 
then we can consider that $f$ has a genuine cusp which is equivalent to $f_{20} f_{02}- f_{11}^2 = 0$ (determinant of Hessian vanishes). Since
 the cusp is genuine cusp so without loss of generality we can assume that $f_{02} \neq 0$. So one can simply 
 construct a curve $\gamma^{i}_{20}(t)  := f + t (\hat{x}-x_i)^2 \eta_i$ for each point $p_i$, where
 $\hat{x}$ is defined below.\\
 Note that 
 \begin{align*}
 \{d\Phi|_{(f, \overline{p})}\} (\gamma^{i^{\prime}}_{{2}, {0}}(0)) = (0, \cdots, \underbrace{*}_\text{$i_{20}$ th position} ,0,0, \cdots, 0)\\
 \end{align*}
 one can observe that this computation proves the claim for cusp.\\
\hf Finally, if $k_i \geq 2$, i.e., $f$ has higher $A_{k_i}$ singularities then we have made a linear change of coordinates so that the kernel of the Hessian is $\partial_x|_{p_i} + m \partial_y|_{p_i}$, 
where $m:= \frac{-f_{xy}(x_i, y_i)}{f_{yy}(x_i, y_i)}$.
Let us now define the curves  
\begin{align*}
\gamma^{i}_{20}(t) & := f + t (\hat{x}-x_i)^2 \eta_i, \hf  \ldots, ~~\gamma^{i}_{k_i 0}(t) := f + t (\hat{x}-x_i)^{k_i} \eta_i,
\end{align*}
for all $i$ from $1$ to $\delta$. Here $\hat{x}:= x + m y$.
 We now note that 
\begin{align*}
\{d\Phi|_{(f, \overline{p})}\} (\gamma^{i^{\prime}}_{\alpha \beta}(0)) 
\end{align*}
span the tangent space of $T_{\mathbf{0}} \mathbb{C}^{c_d}$. This proves the claim. \qed

\begin{lmm}
Restricted to $\H\times L-\mathcal{B}$, the sections  
$\psi_{\textnormal{ev}}$ and $\psi_{\textnormal{T}}$ 
are transverse to zero.
\end{lmm}

\noindent \textbf{Proof:} First, suppose 
\begin{align*}
\psi_{\textnormal{ev}} ([f], \overline{q}, p) &=0 \qquad \iff \qquad f(p) =0. 
\end{align*}
We will produce the following curve. Let us consider a curve $\eta_{00}$ in $\H$ such that $\eta_{00}(p) \neq 0$. 
Consider 
\[ \gamma_{00}(t):= (f + t \eta_{00}, \overline{q}, p). \] 
This proves transversality of the evaluation map. \\ 
\hf \hf Next, let us consider a curve $\eta_{\mathrm{T}}$ such that 
\[ \nabla \eta_{\mathrm{T}}|_p (v) \neq 0 \]
if $v \in T_pL-0$.\\
The construction of a curve $\eta_{\mathrm{T}}$ will follow from above discussion. Now consider the curve 
\[ \gamma_{\mathrm{T}}(t):= f+ t \eta_{\mathrm{T}}.  \]
This proves  transversality of the section $\psi_{\textnormal{T}}$. \qed \\

\noindent Finally, we are ready to prove the main theorem about the multiplicity.

\begin{thm}
\label{mult_th}
Let $\mu \subset \DD$ be the subspace of curves 
passing through $w_d$ generic points and suppose 
\begin{align*}
([f], \overline{q}, p) & \in \H\times L \cap \pi_{\DD}^{-1}(\mu).
\end{align*}
Suppose  
\begin{align}
\psi_{\textnormal{ev}} ([f], \overline{q}, p) &=0, \qquad \psi_{\textnormal{T}} ([f], \overline{q}, p) =0, \qquad 
([f], \overline{q}, p) \in \mathcal{B}^k_{\eta} \cap \Big(\pi_{\D}^{-1}\mu\Big). \label{ord_van}
\end{align}
Then the order of vanishing is $(k+1)$.
\end{thm}

\begin{rem}
Note that $\H\times L \cap \pi_{\DD}^{-1}(\mu)$ is a smooth complex manifold of dimension $2$. Hence it makes sense to talk about the order of vanishing 
of a section of a rank two bundle. 
\end{rem}

\noindent \textbf{Proof:} Suppose $([f], \overline{q}, p)$ satisfies equation \eqref{ord_van}. 
We will construct a neighbourhood of $([f], \overline{q}, p)$ inside $\H \times L$.
Since  
$([f], \overline{q}, p) \in \mathcal{B}^k_{\eta} \cap \Big(\pi_{\D}^{-1}\mu\Big)$ 
and $\mu$ denotes a subspace of curves in $\D$ passing through $w_d$ generic points,   
we conclude that $f$ has an $A_k$ singularity of $p$. Without loss of generality, we can take 
$p:= [0,0,1] \in \mathbb{P}^2$. Let us also assume that the line $L$ passing through $p$ is given by the equation 
\begin{align}\label{Line Equ}
L&:= \{[X,Y,Z] \in \mathbb{P}^2: aX + bY =0\}, 
\end{align}
where $a$ and $b$ are two fixed complex numbers. Let us now write down the Taylor expansion 
of $f$ around the point $p$. Let us define 
\[ x:= \frac{X}{Z} \qquad \textnormal{and} \qquad y:= \frac{Y}{Z}. \] 
Hence, we get that 
\begin{align*}
f &= \frac{f_{20}}{2}x^2 + f_{11} x y + \frac{f_{02}}{2}y^2 + \frac{f_{30}}{6} x^3 + \ldots.  
\end{align*}
If $f$ has an $A_{k \geq 2}$ singularity at $p$, we conclude that $f_{02}$ or $f_{20}$ can not both be zero; 
let us assume in that case $f_{02} \neq 0$. If $f$ has an $A_{1}$ singularity at $p$, then after a linear change of 
coordinates, we can ensure that $f_{02} \neq 0$. Hence, in all the cases, we can assume without loss of generality that $f_{02} \neq 0$. \\
\hf \hf After making a suitable change of coordinates, the function $f$ is given by 
\begin{align*}
f & = \hat{y}^2 + \hat{x}^{k+1}. 
\end{align*}
After the change of coordinates, the line $L$ in \eqref{Line Equ}, will be given by
\begin{align}\label{Final Line}
L&:= \{[X,Y,Z] \in \mathbb{P}^2: \hat{y} + M \hat{x} + E(\hat{y}, \hat{x})=0\},
\end{align}
where $E$ is second order. Note that the line will be either $\hat{y} + M \hat{x} + E(\hat{y}, \hat{x})=0$ or $M^{'} \hat{y} + \hat{x} + E^{'}(\hat{y}, \hat{x})=0$. Without loss of generality, we are assuming that the line is given by \eqref{Final Line}; since $L$ is a generic line, we can assume this (i.e., we are assuming the line is not given by $x =0$). 
Let us now assume that $k$ is even (i.e., $k+1$ odd). A solution to the equation $f=0$,  close to $(0,0)$ is given by 
\begin{align*}
\hat{y}& = t^{k+1}, \qquad \hat{x} = t^2 \qquad \textnormal{$t$ is small but non zero}. 
\end{align*}
Furthermore, every solution to $f =0$ is of this type. We now consider the second equation of evaluating the derivative along $L$. That gives 
us 
\begin{align*}
(M \partial_{\hat{x}} +  \partial_{\hat{y}}) f & =  M f_{\hat{x}} +  f_{\hat{y}} \\ 
                              & = 2 \hat{y} + (k+1) M \hat{x}^{k} \\ 
                              & = 2 t^{k+1} + (k+1) M t^{2 k}.
\end{align*}
Hence, the order of vanishing is $(k+1)$. If $k$ is odd (i.e., $k+1$ is even), then there are two solutions. Each solution vanishes 
with order $\frac{k+1}{2}$; hence the total order of vanishing is $k+1$. In either case, the total order of vanishing is $k+1$. \qed


\section{Explicit Formulas}

\label{ef}
For the convenience of the reader, we will explicitly write down 
the formulas for $N_d^{\mathrm{T}}(A_k)_{1 \leq k \leq 8}$, $N_d^{\mathrm{T}}(A_1 A_k)_{1 \leq k \leq 7}$ and $N_d^{\mathrm{T}}(A_1^{\delta})_{1 \leq \delta \leq 8}$.
These are obtained from Main Theorem (equation \eqref{main_formula}); combined with the numbers given in 
the papers of S.~Basu and R.~Mukherjee (\cite{R.M}, \cite{B.R} and \cite{BM8}). 
We will then use these formulas to make low degree check in section \ref{ldchk}. The formulas for $N_d^{\mathrm{T}}(A_k)_{1 \leq k \leq 8}$ 
are:
\begin{align*}
N_d^{\mathrm{T}}(A_1)& =6d(d-1)(d-2), \qquad N_d^{\mathrm{T}}(A_2) =12 ( 2 d^3- 8 d^2+ 8 d -1), \\ 
N_d^{\mathrm{T}}(A_3)& = 4 (25 d^3 - 146 d^2+ 228 d- 84), \qquad  N_d^{\mathrm{T}}(A_4) = 120 (3 d^3- 20 d^2+ 36 d -15), \\ 
N_d^{\mathrm{T}}(A_5)&= 36 (35 d^3- 260 d^2+ 524 d  -239),  \qquad N_d^{\mathrm{T}}(A_6)= 7 (632 d^3 - 5134 d^2+ 11343 d -5538), \\ 
N_d^{\mathrm{T}}(A_7)&= 24 (651 d^3 - 5702 d^2+ 13602 d -7002) \qquad \textnormal{and} \\ 
N_d^{\mathrm{T}}(A_8)&= 288 (190 d^3- 1778 d^2+ 4533 d -2436).
\end{align*}
Next, the formulas for $N_d^{\mathrm{T}}(A_1 A_k)_{1 \leq k \leq 7}$ are: 
\begin{align*}
N_d^{\mathrm{T}}(A_1^2)& = 2(9 d^5- 45 d^4+ 30 d^3+ 123 d^2- 145 d +6), \\ 
N_d^{\mathrm{T}}(A_1 A_2) & = 12 (d-3) (6 d^4-18 d^3-22 d^2+67 d-13), \\ 
N_d^{\mathrm{T}}(A_1 A_3) & = 12 (25 d^5-171 d^4+187 d^3+774 d^2-1535 d+426), \\ 
N_d^{\mathrm{T}}(A_1 A_4) & = 20 (54 d^5-414 d^4+534 d^3+2238 d^2-5207 d+1815), \\ 
N_d^{\mathrm{T}}(A_1 A_5) & = 18 (210 d^5-1770 d^4+2572 d^3+11299 d^2-29650 d+11959), \\ 
N_d^{\mathrm{T}}(A_1 A_6) & = 21 (632 d^5-5766 d^4+9164 d^3+42837 d^2-123391 d+55068), \\ 
N_d^{\mathrm{T}}(A_1 A_7) & = 8 (5859 d^5-57177 d^4+97677 d^3+485874 d^2-1509623 d+725940).
\end{align*}
Finally, the formulas for $N_d^{\mathrm{T}}(A_1^{\delta})_{3 \leq \delta \leq 8}$ are: 
\begin{align*}
N_d^{\mathrm{T}}(A_1^3) & = 6 (9 d^7 -63 d^6+36 d^5+549 d^4-857d^3-1148 d^2 + 2266 d -300), \\ 
N_d^{\mathrm{T}}(A_1^4) & = 18 (9 d^9-81 d^8+36 d^7+1458 d^6-2834 d^5-8500 d^4+22455 d^3+13543 d^2-49222 d+10488), \\ 
N_d^{\mathrm{T}}(A_1^5) & = 6 \big(81 d^{11}-891 d^{10}+ 270 d^9 + 27270 d^8- 63450 d^7-303912 d^6\\ 
 & ~~ +1014807 d^5+ 1348725 d^4- 6097876 d^3- 1168832 d^2   + 12259248 d -3513840\big), \\ 
N_d^{\mathrm{T}}(A_1^6) & = 1458 d^{13} - 18954 d^{12} + 2916 d^{11}+ 882090 d^{10}- 2390310 d^9 - 
 15901596 d^8+ 64328418 d^7 \\ 
 & + 130916898 d^6- 732619008 d^5-395637750 d^4  + 3855455766 d^3\\ 
 & - 418407408 d^2 - 7418026440 d + 2643818400, \\
N_d^{\mathrm{T}}(A_1^7) & = 4374 d^{15}- 65610 d^{14}+ 4317138 d^{12}- 13352850 d^{11}  - 
 114293592 d^{10}+ 543520530 d^9\\ 
 &+ 1481762970 d^8 - 9946281060 d^7 - 
 8470208502 d^6 + 95900422338 d^5+ 1014814332 d^4\\  
 & -467415101124 d^3 + 168796887984 d^2 + 880782565392 d -374053619520 \qquad \textnormal{and} \\ 
N_d^{\mathrm{T}}(A_1^8) & = 13122 d^{17}- 
 223074 d^{16}- 34992 d^{15} + 19717992 d^{14}- 68543496 d^{13} \\ 
 & -719400528 d^{12} + 3933317556 d^{11}+ 13400193204 d^{10}- 
 105120249336 d^9 \\ 
 & - 119845037160 d^8+ 1587321808632 d^7             + 
 150918108768 d^6 - 13835625254910 d^5 \\ 
 & + 5746599271062 d^4 + 
 64281794069664 d^3 - 38151916883064 d^2 \\ 
 & - 120388035085920 d+59358641529600.
\end{align*}

\section{Low degree checks} 

\label{ldchk}
In this section, we will make  some non trivial low degree checks by comparing our formulas with the results of others.

\subsection{Verification with the Caporaso-Harris formula}  
We will start by verifying the numbers $N_d(A_1^{\delta})_{1 \leq \delta \leq 8}$. We note that the Caporaso-Harris formula (obtained in \cite{CH}) 
computes $N_d(A_1^{\delta})^{\mathrm{T}}$ for any $\delta$. 
We have verified that our formulas for $N_d(A_1^{\delta})_{1 \leq \delta \leq 8}$ produce the same answer as the 
Caporaso-Harris formula for several values of $d$; we have written a C++  program to implement the Caporaso-Harris formula (which is available 
on request). The reader is invited to use the C++ program to check that it produces the same answer given by our formula (explicitly written down 
in section \ref{ef}) for any specific value of $d$.

\subsection{Verification of $N^{\mathrm{T}}_d(A_2)$ using a result of Kazaryan}
In \cite{Kaz}, Kazaryan has computed the number $N_d(A_1 A_2 A_3)$, the characteristic number of degree $d$ curves 
with one node, one cusp and one tacnode. According to Kazaryan's formula, that number is $2256$ when $d=4$. 
We will verify that number.\\   
\hf \hf We note that $N_4(A_1 A_2 A_3)$ is the number of quartics through $8$ points that have one node, one cusp and one tacnode. This 
can happen if the curve breaks into a cubic and a line, such that the cubic has a cusp and is tangent to the given line (and the entire 
configuration passes through $8$ points). Since the 
cubic is tangent to the given line, it will intersect the curve at one more  point. Let us now find out how many such configurations are there. 
First of all, we could place a line through $2$ points and a cuspidal cubic through $6$ points tangent to a given line. There are a total of 
\[ \binom{8}{2} \times N^{\mathrm{T}}_{3}(A_2) \] 
such configurations. The other possibility is that we place a cuspidal cubic through $7$ points and a line through one point that 
is tangent to this cuspidal cubic. We claim that the total number of such configurations ($n$) is 
\[ n= 3 N_3(A_2). \] 
We will justify this shortly. 
Using the values of $N^{\mathrm{T}}_{3}(A_2)$ and $N_{3}(A_2)$, we note that 
\begin{align*}
\binom{8}{2} \times N^{\mathrm{T}}_{3}(A_2) + \binom{8}{7} \times n & = 
\binom{8}{2} \times N^{\mathrm{T}}_{3}(A_2) + \binom{8}{7} \times 3N_3(A_2)  \\ 
& = \binom{8}{2} \times 60 + \binom{8}{7} \times 3\times 24 \\
& = 2256.
\end{align*}
This agrees with the number predicted by Kazaryan's formula. \\  
\hf \hf Let us now justify the value of $n$. 
Let us denote $\DD_1$ and $\DD_3$ to be the space of lines and space of cubics in $\mathbb{P}^2$ respectively. 
We note that $\DD_1$ and $\DD_3$ are isomorphic to $\mathbb{P}^2$ and $\mathbb{P}^9$ respectively. \\
  
\hf \hf Let us define 
\begin{align*}
\mathcal{S} &:= \{([f], q) \in \DD_3 \times \mathbb{P}^2: ~\textnormal{$f$ has an $A_2$ singularity at $q$}\}. 
\end{align*}
For notational convenience, let us denote $\mathbb{P}^2_1$ and $\mathbb{P}^2_3$ to be two isomorphic copies of $\mathbb{P}^2$.  
With that notation, we define the following space 
\begin{align}
\mathcal{Z} &:= \{ ([f_1], q_1, [f_3], q_3) \in \DD_1 \times \mathbb{P}^2_1 \times \DD_3 \times \mathbb{P}^2_3: ~~([f_3], q_3) \in \overline{\mathcal{S}}, 
~~f_1(q_1) =0, ~~f_3(q_1) =0\}. \label{z_def} 
\end{align}
Next, we note that over the space $\mathcal{Z}$, we have the following short exact sequence of bundles  
\begin{equation}
\label{ses_line_bundle}
\xymatrix{
0 \ar[r] & \mathbb{L}:= \textnormal{Ker}(\nabla f_1|_{q_1}) \ar[r] & T\P^2|_{q_1} \ar[r]^-{\nabla f_1|_{q_1}} & \gamma_{\mathcal{D}_1}^* \otimes \gamma_{\P^2}^* \lra 0. 
} 
\end{equation}
Let us now define the following set 
\begin{align}
X &:= \{([f_1],q_1, [f_3], q_3) \in \DD_1 \times \mathbb{P}^2_1 \times \DD_3  \times \mathbb{P}^2_3: ([f_1],q_1, [f_3], q_3) \in \mathcal{Z},
~~\nabla f_3|_{q_{1}} (v) = 0, ~~\forall v \in \mathbb{L}\}. \label{x_def}
\end{align}
Let us now denote $y_1$, $y_3$, $a_1$ and $a_3$ to be the hyperplane classes of $\DD_1$, $\DD_3$, $\mathbb{P}^2_1$ and $\mathbb{P}^2_3$ respectively.  
We note that intersecting $[X]$ with $y_3$ corresponds to studying the subspace of cubics passing through a generic point 
and  intersecting $[X]$ with $y_1$ corresponds to studying the subspace of lines passing through a generic point. 
Our aim is to count the configurations where the cubic passes through $7$ points and the line passes through $1$ point. 
Hence, let us intersect $[X]$ this with $y_1 y_3^7$. However, this intersection will also include the number of lines that 
pass through the given point and the cuspidal point of the cubic. By using the same argument as in the proof of Theorem \ref{mult_th}, 
this configuration contributes with a multiplicity of $3$. Hence, 
\begin{align}
[X]\cdot [y_1 y_3^{7}] & = n + 3N_{3}(A_2). \label{n_def}
\end{align}
It remains to compute $[X]\cdot [y_1 y_3^{7}]$. 
First we note that 
\begin{align}
\langle y_3^7, ~[\overline{\mathcal{S}}]\rangle = N_3(A_2).\label{y17} 
\end{align}
Using equations \eqref{z_def}, \eqref{x_def}, \eqref{ses_line_bundle} and \eqref{y17}, we conclude that 
\begin{align}
[X]\cdot [y_1 y_3^{7}] & = \langle y_1 y_3^7 e(\gamma_{\mathcal{D}_1}^* \otimes \gamma_{\P^2_1}^*)e(\gamma_{\mathcal{D}_3}^* 
\otimes \gamma_{\P^2_3}^{*3}) e(\gamma_{\mathcal{D}_3}^* \otimes \gamma_{\P^2_3}^{*3} 
\otimes \mathbb{L}^*), ~\DD_1 \times \mathbb{P}^2_1 \times \overline{\mathcal{S}} \rangle \nonumber \\
                              & = 6 N_3(A_2). \label{x_eul}
\end{align}
Equations \eqref{n_def} and \eqref{x_eul}, we conclude that $n = 3 N_3(A_2)$ as claimed.

\begin{rem}
It should be possible to generalize our method to get a similar type of result for certain other types of 
singularities such as $D_k$, $E_6$, $E_7$, and $E_8$. We are not aware of any low degree checks involving tangency conditions 
with $D_k$, $E_6$, $E_7$, and $E_8$ singularities which will 
support any prediction.
We also hope that this method can be employed to generalize these results for other complex surfaces.
\end{rem}

\section{Acknowledgement}
I am indebted to my advisor Ritwik Mukherjee for introducing me to the subject and several fruitful discussions indicating how this problem related to 
several other interesting as well as important branches in mathematics. I want to thank Somnath Basu for several fruitful discussions. I also thank Nilkantha Das and Soumya Pal for writing a C++ and Python program to implement Caporaso-Harris 
formula for verification. This work was supported by the Council of Scientific and Industrial Research (CSIR award Sr.No. 1121540885).

\address{School of Mathematics, National Institute of Science Education and Research, Bhubaneswar, HBNI, Odisha 752050, India.}
 	\email{anantadulal.paul@niser.ac.in } \\
 	\address{Homi Bhaba National Institute, Training School Complex, Anushakti Nagar, Mumbai, 400094, India.}
\end{document}